\documentclass[a4paper,10pt]{article}

\usepackage{amssymb,amsmath,amscd,latexsym,url}
\usepackage{setspace}
\newtheorem{defi}{Definition}[section]
\newtheorem{theo}[defi]{Theorem}
\newtheorem{lem}[defi]{Lemma}

\newtheorem{coro}[defi]{Corollary}
\newtheorem{rema}[defi]{Remark}

\def\RR{\mathbb{R}}
\def\NN{\mathbb{N}}
\def\Ccal{\mathcal{C}}

\def\adom{\overline{\Omega}}
\def\bom{\partial \Omega}

\def\proof{\emph{Proof :} \hspace*{0.5cm}}
\newcommand{\fin}{\mbox{}\nolinebreak\hfill\rule{2mm}{2mm}\par\medbreak}

\newcommand{\refeq}[1]{$(\ref{#1})$}

\newcommand{\norme}[1]{\parallel #1 \parallel _\infty}
\def\dt{\partial_t}
\def\dn{\partial_\nu}

\title{The Fujita phenomenon in exterior domains under dynamical boundary conditions II\footnote{We improve the previous results from The Fujita phenomenon in exterior domain under dynamical boundary conditions, \emph{Technical Report L.M.P.A.} \textbf{383}, October 2008.} \footnote{To appear in Asymptotic Analysis.}} 
\author{Jean-Fran\c cois Rault}
\date{}

\begin{document}

\bibliographystyle{amsplain}

\maketitle

\begin{center}

LMPA Joseph Liouville, FR 2956 CNRS, \\
Universit\'e du Littoral C\^ote d'Opale \\
50 rue F. Buisson, B.P. 699, F-62228 Calais Cedex (France)\\
\url{jfrault@lmpa.univ-littoral.fr}\\
\end{center}
\bigskip

\begin{abstract}
\noindent The Fujita phenomenon for nonlinear parabolic problems $\partial_t u=\Delta u +u^p$ in an exterior domain of $\RR^N$ under dissipative dynamical boundary conditions $ \sigma \partial_t u+\partial_\nu u=0$ is investigated in the superlinear case.  As in the case of Dirichlet boundary conditions (see Refs. \cite{BL1} and \cite{MS}), it turns out that there exists a critical exponent $p=1+\frac{2}{N}$ such that blow-up of positive solutions always occurs for subcritical exponents, whereas in the supercritical case global existence can occur for small non-negative initial data.\\


\end{abstract}
\bigskip

\section{Introduction}

Let $\Omega$ be an exterior domain of $\RR^N$, that is to say a connected open set $\Omega$ such that $\adom^c$ is a bounded domain when $N\geq 2$, and in dimension one, $\Omega$ is the complement of a real closed interval. We always suppose that the boundary $\bom$ is of class $\Ccal^2$. The outer normal unit vector field is denoted by $\nu: \bom \rightarrow \RR^N$ and the outer normal derivative by $\dn$. Let $p$ be a real number with $p>1$ and $\varphi$ be a continuous function in $\adom$. Consider the following nonlinear parabolic problem
\begin{equation}\label{Ext_pbm}
\left\{
    \begin{array}{ll}
        \partial _t u = \Delta u + u^p& \textrm{ in $\adom \times (0,\infty)$}, \\
        \mathcal{B} _\sigma (u) := \sigma \partial _t u + \partial _\nu u= 0 & \textrm{ on $\bom \times (0,\infty)$}, \\
        u(\cdot , 0) = \varphi & \textrm{ in $\adom$ }.
    \end{array}
\right.
\end{equation}
The aim of this paper is to show that the well-known Fujita phenomenon in the case of $\Omega=\RR^N$ (see Ref. \cite{Fujita}) and in the case of Dirichlet boundary conditions (see Refs. \cite{BL1} and \cite{MS}) still holds for the dynamical boundary conditions. One can notice that dynamical boundary conditions $\mathcal{B} _\sigma(u) =0$ with $\sigma \equiv 0$ correspond to the Neumann boundary conditions, which case has been discussed by Levine and Zhang \cite{LZ}. It is already known, by Bandle, von Below and Reichel in \cite{BBR1}, that for $p \in (1, 1+\frac{2}{N})$, also for $p=1+\frac{2}{N}$ if $N\geq 3$, and for constant coefficient $\sigma \in [0,\infty)$, all positive solutions of \refeq{Ext_pbm} blow up in finite time. In addition, if the complement is star-shaped there exist global positive solutions of class $\Ccal^1$ for $p>1+\frac{2}{N}$  by \cite{BBR1}. Our purpose is to show the existence of global positive solutions of Problem \refeq{Ext_pbm} for sufficiently small initial data in the supercritical case ($p>1+\frac{2}{N}$) for any exterior domain. Moreover our condition on $\sigma$ is more general. Throughout, we shall assume the dissipativity condition
\begin{equation}\label{dissip}
\sigma \geq 0 \textrm{ on } \bom \times (0,\infty)
\end{equation}
and dealing with classical solutions
\begin{equation}\label{sreg}
\sigma \in \Ccal^1( \bom \times (0,\infty)).
\end{equation}
The initial data is always supposed to be continuous, non-trivial, bounded, non-negative in $\adom$,  and vanishing at infinity:
\begin{equation}\label{i-data}
\varphi \in \Ccal(\adom), \ 0<\norme{\varphi} < \infty, \ \varphi \geq 0,\ \lim_{\parallel x \parallel_2 \rightarrow\infty}\varphi (x)=0.
\end{equation}
In the case $\Omega = \RR^N$, the boundary condition is dropped and the result is well known by the classical paper of Fujita \cite{Fujita}. Thus, we will suppose $\Omega \not= \RR^N$.

\section{Preliminaries}\label{existence}
First, we give the definition of positive solution which is understood along this paper.
\begin{defi}
A positive solution of Problem \refeq{Ext_pbm} is a positive function $u: (x,t) \mapsto u(x,t)$ of class $\Ccal(\adom \times [0,T)) \cap \Ccal^{2,1}(\adom \times (0,T))$, satisfying
\begin{displaymath}
\left\{
    \begin{array}{ll}
        \partial _t u = \Delta u + u^p& \textrm{ in $\adom \times (0,T)$}, \\
        \mathcal{B} _\sigma (u) := \sigma \partial _t u + \partial _\nu u= 0 & \textrm{ on $\bom \times (0,T)$}, \\
        u(\cdot , 0) = \varphi & \textrm{ in $\adom$ },
    \end{array}
\right.
\end{displaymath}
where $\varphi$ is a function, given in $\Ccal(\adom)$.  
The time $T  \in [0,\infty]$ is the maximal existence time of the solution $u$. If $T=\infty$, the solution u is called global.
\end{defi}
From \cite{BL1}, if $T<\infty$, $u$ blows up in finite time, that is to say:
\begin{displaymath}
\lim_{t\nearrow T} \sup_{x\in \adom} u(x,t)= \infty.
\end{displaymath}
Note that for initial data $\varphi$ of class $\Ccal^2(\adom)$, the solution $u$ is $\Ccal^{2,1}(\adom \times [0,T))$, whereas $u \in \Ccal(\adom \times [0,T)) \cap \Ccal^{2,1}(\adom \times (0,T))$ if $\varphi$ is only continuous in $\adom$. Then, let us recall a standard procedure to construct solutions of Problem \refeq{Ext_pbm} in outer domains for uniformly bounded and continuous initial data $\varphi$. Let $B(0,R)$ be the ball centered at the origin of radius $R>0$ such that $\adom^c \subset B(0,R)$. For any $n \in \NN$, we set $B_n :=B(0,R+n)$ and  $\Omega_n := \Omega \cap B_n$. The boundary of $\Omega_n$ is decomposed into two disjoint open sets:
\begin{displaymath}
\bom _n = \bom \dot{\cup} \partial B_n.
\end{displaymath}
Define also an increasing sequence of initial data $(\varphi_n)_{n\in \NN^*}$ such that
\begin{eqnarray}\label{init-seq}
0\leq \varphi_n \leq \varphi & \textrm{ in }& \adom_n , {} \nonumber \\
 \varphi_n \equiv 0  &\textrm{ on }& \partial B_n, {}  \\
 \varphi_n = \varphi    & \textrm{ in } &\adom_{n-1},{}  \nonumber
\end{eqnarray} 
and consider the following problem with mixed boundary conditions
\begin{displaymath}\tag{$P(n)$}
\left\{
    \begin{array}{ll}
        \partial _t u = \Delta u + u^p& \textrm{ in $\overline{\Omega}_n \times (0,\infty)$}, \\
        \mathcal{B} _\sigma (u) := \sigma \partial _t u + \partial _\nu u= 0 & \textrm{ on $\partial \Omega \times (0,\infty)$}, \\
        u= 0 & \textrm{ on $\partial B_n \times (0,\infty)$}, \\
        u(\cdot , 0) = \varphi_n & \textrm{ in $\overline{\Omega}_n$ }.
    \end{array}
\right.
\end{displaymath}
Let $z$ be the maximal solution of
\begin{displaymath}
\left\{
    \begin{array}{ll}
        \dot{z} = z^p, \\
        z(0) = \norme{\varphi},
    \end{array}
\right.
\end{displaymath}
with maximal existence time $t_0=\frac{1}{(p-1) \norme{\varphi}^{p-1}}$. It is known from \cite{BP1} that, for each $n \in \NN^*$, Problem $(P(n))$ has a solution $u_n \in \Ccal(\adom_n \times [0,T_n)) \cap \Ccal^{2,1}(\adom_n \times (0,T_n))$, where $T_n$ is the maximal existence time of $u_n$. Moreover by comparison principle from \cite{BDC1}, we have, for any $n\in \NN^*$, $0\leq u_n \leq u_{n+1}$ and $u_n(\cdot,t) \leq z(t)$ in $\adom_n$, so we have also $t_0\leq T_n$. Hence we obtain a sequence $(u_n)_{n\in \NN^*}$ of functions in $\Ccal(\adom_n \times [0,t_0)) \cap \Ccal^{2,1}(\adom_n \times (0,t_0))$. Then, standard arguments based on a priori estimates for the heat equation imply  $u_n \rightarrow u$ in the sense of $\Ccal^{2,1}_{loc}(\adom \times (0,t_0))$ as $n\rightarrow \infty$, where $u$ is a positive solution of Problem \refeq{Ext_pbm}, see Refs. \cite{BBR1} and \cite{LSU}. Moreover, since $u_n$ vanishes on $\partial B_n$ for each $n\in \NN^*$, the solution $u$ vanishes at infinity: 
\begin{displaymath}
\lim_{\parallel x \parallel_2 \rightarrow \infty} u(x,t) = 0 \ , \forall \ t \in (0,T) \ .
\end{displaymath}
\noindent Note that $t_0$ is only a lower bound  for the maximal existence time of solutions $u_n$ and $u$, and it is possible that the times $T_n$ and $T$ are infinite. Indeed, results on blow-up for problems under dynamical boundary conditions from \cite{BP2} can not be applied to the problems $(P(n))$ with mixed boundary conditions because their solutions $u_n$ vanish on a part of the boundary.

\section{Global existence in dimension $N \geq 3$}

\noindent Throughout this section, we consider supercritical exponent $p$:
\begin{displaymath}
p>1+\frac{2}{N}.
\end{displaymath}
Our technique will be to construct a function that bounds from above each solution $u_n$ of Problem $(P(n))$. This will give us a sequence $(u_n)_{n\in \NN^*}$ of global solutions in $\Ccal(\adom_n \times [0,\infty)) \cap \Ccal^{2,1}(\adom_n \times (0,\infty))$, thus the solution $u$ of \refeq{Ext_pbm} must be global too. We will proceed by using the solution of the Neumann Problem 
\begin{equation}\label{Ext_Neu}
\left\{
    \begin{array}{ll}
        \partial _t v = \Delta v + v^p& \textrm{ in $\adom \times (0,\infty)$}, \\
        \partial _\nu v= 0 & \textrm{ on $\bom \times (0,\infty)$}, \\
        v(\cdot , 0) = \psi & \textrm{ in $\adom$ },
    \end{array}
\right.
\end{equation}
with $\psi$ verifying \refeq{i-data}. In \cite{LZ}, Levine and Zhang proved that Problem \refeq{Ext_Neu} admits global positive solutions for sufficiently small initial data. We show that the solution $v$ of Problem \refeq{Ext_Neu} bounds from above the solution $u$ of Problem \refeq{Ext_pbm} if the initial data are well ordered ($\varphi \leq \psi$) and if $\psi$ satisfies the following hypotheses:
\begin{equation}\label{idN1}
\psi \in \Ccal^2(\adom) \ ,
\end{equation}
and for every $n\in \NN^*$
\begin{equation}\label{idN2}
\Delta \psi_n + \psi_n^p \geq 0 \textrm{ in } \adom_n \ ,
\end{equation}
where $(\psi_n)_{n \in \NN^*}$ is the sequence of truncated  initial data introduced in \refeq{init-seq}. We need a technical lemma, similar to Lemma 2.1 of \cite{BP1}.

\begin{lem}
Let $\psi$ be a function satisfying \refeq{i-data}, \refeq{idN1} and \refeq{idN2}. For every $n\in \NN^*$, the solution $v_n$ of Problem $(P(n))$ under the Neumann boundary conditions and with the truncated initial data $\psi_n$ verifies:
\begin{displaymath}
\dt v_n \geq 0 \textrm{ in } \adom_n \times (0,T_{v_n}),
\end{displaymath}
where $T_{v_n}$ is the maximal existence time of $v_n$.
\end{lem}
\proof The function $v_n$ is solution of the following problem
\begin{displaymath}
\left\{
    \begin{array}{ll}
        \partial _t v_n = \Delta v_n + v_n^p& \textrm{ in $\adom_n \times (0,T_{v_n})$}, \\
        \partial _\nu v_n= 0 & \textrm{ on $\bom \times (0,T_{v_n})$}, \\
         v_n= 0 & \textrm{ on $\partial B_n \times (0,T_{v_n})$}, \\
        v_n(\cdot , 0) = \psi_n & \textrm{ in $\adom_n$ },
    \end{array}
\right.
\end{displaymath}
in the bounded domain $\adom_n$. From \refeq{i-data} and from the strong maximum principle in \cite{BDC1}, we claim 
\begin{displaymath}
 v_n > 0 \textrm{ on } (\Omega_n \cup \bom) \times (0,T_{v_n}).
\end{displaymath}
Then, by regularity results from \cite{LSU}, we obtain $v_n \in \Ccal^{2,2}(\adom_n \times (0,T_{v_n})$, and for $y= \dt v_n \in \Ccal^{2,1}(\adom_n \times (0,T_{v_n})$ we have
\begin{displaymath}
\left\{
    \begin{array}{ll}
        \partial _t y = \Delta y + p v_n^{p-1}y& \textrm{ in $\adom_n \times (0,T_{v_n})$}, \\
        \partial _\nu y= 0 & \textrm{ on $\bom \times (0,T_{v_n})$}, \\
         y= 0 & \textrm{ on $\partial B_n \times (0,T_{v_n})$},
    \end{array}
\right.
\end{displaymath}
and  $y(\cdot , 0) \geq 0$ in $\adom_n$ thanks to \refeq{idN2}. By the comparison principle in \cite{BDC1}, we conclude: $y \geq 0$ in $\adom_n \times (0,T_{v_n})$. \\
\fin

\begin{lem}\label{u < Neumann}Let a coefficient $\sigma$ verifying \refeq{dissip} and \refeq{sreg}, two functions $\varphi$ and $\psi$ satisfying \refeq{i-data} and $\psi$ with \refeq{idN1} and \refeq{idN2}. If 
\begin{equation}\label{comparaison}
\varphi \leq \psi \textrm{ in } \adom,
\end{equation}
then  Problem \refeq{Ext_pbm} with initial data $\varphi$ admits a solution $u$ verifying
\begin{displaymath}
u \leq v \textrm{ in } \adom \times (0,T_v),
\end{displaymath}
and 
\begin{displaymath}
0<T_v \leq T \leq \infty,
\end{displaymath}
where $v$ is solution of Problem \refeq{Ext_Neu} with initial data  $\psi$, of maximal existence time $T_v$.
\end{lem}
\proof We consider the sequences of truncated solutions $(u_n)_{n\in\NN^*}$ and $(v_n)_{n\in\NN^*}$ respectively associated to the solutions $u$ and $v$. Let $n \in \NN^*$. First, we show that $v_n \leq v$ in $\adom_n \times [0,T_v)$. By construction \refeq{init-seq}, we have $\psi_n \leq \psi$ in $\adom_n$. Since $v$ is a positive solution of Problem \refeq{Ext_Neu}, it satisfies
\begin{displaymath}
\left\{
    \begin{array}{ll}
        \partial _t v \geq \Delta v + v^p& \textrm{ in $\adom_n \times (0,T_v)$}, \\
        \dn v \geq 0 & \textrm{ on $\bom \times (0,T_v)$}, \\
         v \geq 0 & \textrm{ on $\partial B_n \times (0,T_v)$}, \\
        v(\cdot , 0) \geq \psi_n & \textrm{ in $\adom_n$ }.
    \end{array}
\right.
\end{displaymath}
As $v_n$ is a positive solution of $(P(n))$ under Neumann boundary conditions, we obtain from the comparison principle in \cite{BDC1}
\begin{equation}\label{VnV}
v_n \leq v \textrm{ in }\adom_n \times [0,\tau),
\end{equation} 
for all $0< \tau < \min \{ \ T_{v_n}, \ T_v \ \}$. We deduce $T_v\leq T_{v_n}$. Then, we show that $u_n \leq v_n$ in $\adom_n \times [0,T_{v_n})$. The previous lemma ensures that $\dt v_n \geq 0$. From \refeq{dissip}, we obtain
\begin{displaymath}
\sigma \dt v_n +\dn v_n \geq 0 \textrm{ on } \bom \times (0,T_{v_n}).
\end{displaymath}
Next, $v_n$ is a positive solution of $(P(n))$, and thanks to \refeq{comparaison}, $v_n$ verifies
\begin{displaymath}
\left\{
    \begin{array}{ll}
        \partial _t v_n \geq \Delta v_n + v_n^p& \textrm{ in $\adom_n \times (0,T_{v_n})$}, \\
        \sigma \dt v_n +\dn v_n \geq 0 & \textrm{ on $\bom \times (0,T_{v_n})$}, \\
         v_n \geq 0 & \textrm{ on $\partial B_n \times (0,T_{v_n})$}, \\
        v_n(\cdot , 0) \geq \varphi_n & \textrm{ in $\adom_n$ }.
    \end{array}
\right.
\end{displaymath}
Again, by the comparison principle in \cite{BDC1} and by definition of $u_n$, we obtain
\begin{equation}\label{unVn}
u_n \leq v_n \textrm{ in }\adom_n \times [0,\tau),
\end{equation} 
for all $0<\tau<\min \{ \ T_n, \ T_{v_n} \ \}$, and hence $T_{v_n} \leq T_{n}$. From equations \refeq{VnV} and \refeq{unVn}, we have $T_v\leq T_n$ and $u_n \leq v$ in $\adom_n \times [0,T_v)$. Thus the solution $u$ of Problem \refeq{Ext_pbm}, obtained as the limit of the sequence $(u_n)_{n\in \NN^*}$ with the procedure described in section \ref{existence}, verifies $u \leq v$ in $\adom \times [0,T_v)$ and $T_{v} \leq T$.\\
\fin

\noindent 
\noindent Now, we just have to choose an initial data $\psi_*$ sufficiently small such that Problem \refeq{Ext_Neu} admits a global positive solution (see Ref. \cite{LZ}), and satisfying \refeq{idN1} and \refeq{idN2}.

\begin{theo} Under conditions \refeq{dissip}, \refeq{sreg} and \refeq{i-data}, Problem \refeq{Ext_pbm} admits global positive solutions for sufficiently small initial data. Moreover, some of these solutions vanish at infinity.
\end{theo}
\proof An initial data $\varphi$ verifying \refeq{i-data} and with $\varphi \leq \psi_*$ in $\adom$ allows us to conclude thanks to Lemma \ref{u < Neumann}.\\ 
\fin

\begin{rema} One can notice that only the dissipativity and the regularity of the coefficient $\sigma$ are needed. We are not obliged to impose any restriction like $\sigma$ bounded or $\dt \sigma \equiv 0$. Moreover, the hypotheses \refeq{idN1} and \refeq{idN2} on the initial data  $\psi$ of Problem \refeq{Ext_Neu} are strictly technical and do not concern the initial data $\varphi$ of Problem \refeq{Ext_pbm}. 
\end{rema}

\section{Global existence in lower dimension}

\noindent In this case, we can not use Levine and Zhang's result because it is proved only for dimension $N\geq 3$: they used some estimates for Green's functions, specific to dimension $N\geq 3$. We need an additional hypothesis on the coefficient $\sigma$. There exists a constant $\varsigma \in [0,\infty)$ such that
\begin{equation}\label{sborne}
\forall (x,t) \in \bom \times [0,\infty) \ : \ \sigma (x,t) \leq \varsigma \ .
\end{equation}
We begin with the  case of dimension $2$. Until now, Bandle - von Below - Reichel's lemma, concerning star-shaped domains, is the best result:
\begin{lem} \emph{\cite{BBR1}, Lemma 28. } \label{dom_etoile}
Suppose $\sigma$ is a positive constant. If $\Omega^C$ is star-shaped with respect to the origin and if $\displaystyle \min_{\partial \Omega} | x \cdot \nu(x) | \geq \sigma N$, then there exist positive global solutions of Problem \refeq{Ext_pbm}, which vanish at infinity, for sufficiently small initial data. 
\end{lem}

\noindent This allows us to deduce the following result for problems with mixed boundary conditions.

\begin{coro}\label{sol_mixte}
Suppose conditions \refeq{dissip}, \refeq{sreg} and \refeq{sborne}. Let $y \in \bom$. There exists a neighborhood $N_y$ of $y$ relatively open in $\bom$ such that  the following parabolic problem with mixed boundary conditions 
\begin{displaymath}
\left\{
    \begin{array}{ll}
        \partial _t u = \Delta u + u^p& \textrm{ in $\overline{\Omega} \times (0,\infty)$}, \\
        \mathcal{B} _\sigma (u) = 0 & \textrm{ on $N_y \times (0,\infty)$}, \\
        u= 0 & \textrm{ on $\bom \setminus N_y \times (0,\infty)$}, \\
        u(\cdot , 0) = \varphi & \textrm{ in $\overline{\Omega}$ }.
    \end{array}
\right.
\end{displaymath}
admits global positive solutions which vanish at infinity, for sufficiently small initial data $\varphi$ satisfying \refeq{i-data}.
\end{coro}
\proof Let $\mu$ be a vector in $\RR^N$ such that the scalar product between the vector $(y+\mu)$ and the outer normal unit vector at $y$ satisfies
\begin{equation}\label{mu_bom}
(y + \mu)\cdot \nu (y) < - \varsigma N.
\end{equation}
Then, as the mapping ($\bom \ni z \mapsto (z+\mu) \cdot \nu(z) \in \RR$) is continuous, the above inequality remains true on an open neighborhood $N_y \subseteq \bom$ of $y$. We obtain the statement of the corollary by using the comparison principle from \cite{BDC1} and the function $U$ defined on $\adom \times [0,\infty)$ by
\begin{displaymath}
U(x,t)= A (t+1)^{-\gamma} \exp{\frac{-\parallel x + \mu \parallel_2^2}{4(t+1)}}, 
\end{displaymath}
with $\displaystyle A= \frac{1}{2}\Big(\frac{N}{2} -\frac{1}{p-1} \Big)^\frac{1}{p-1}$ and $\gamma = \frac{1}{p-1}$. It is clear that $U \geq 0$, belongs to $\Ccal^{2,1}(\adom \times [0,\infty)) $ and satisfies:
\begin{eqnarray}
\dt U    & = & \Big( \frac{-\gamma}{t+1} + \frac{\parallel x +\mu \parallel_2^2}{4(t+1)^2} \Big) U, {} \nonumber \\
\Delta U & = & \Big( \frac{-N}{2(t+1)} + \frac{\parallel x + \mu \parallel_2^2}{4(t+1)^2} \Big) U,{}  \nonumber\\
\dn U    & = & \Big( \frac{-(x+\mu)\cdot\nu(x)}{2(t+1)} \Big) U.{}  \nonumber
\end{eqnarray} 
We have:
\begin{displaymath}
\dt U - \Delta U = \Big( \frac{-\gamma}{t+1} + \frac{N}{2(t+1)} \Big) U = \Big( \frac{-2 \gamma +N}{2(t+1)} \Big) U\ ,
\end{displaymath}
and by definition of the constants $A$ and $\gamma$, we obtain
\begin{displaymath}
\dt U - \Delta U -U^p \geq 0 \textrm{ in } \adom \times [0, \infty) \ .
\end{displaymath}
Then, on $\bom$, we have 
\begin{displaymath}
\sigma \dt U +\dn U = \Big( \frac{-2\sigma \gamma - (x+\mu)\cdot\nu(x) }{2(t+1)} + \frac{\sigma \parallel x +\mu \parallel_2^2}{4(t+1)^2} \Big) U.
\end{displaymath}
Since $p>1+\frac{2}{N}$, using \refeq{sborne} and ignoring the non-negative term $\frac{\sigma \parallel x+\mu \parallel_2^2}{4(t+1)^2}$, we can reduce the last equation to:
\begin{displaymath}
\sigma \dt U +\dn U \geq \Big( \frac{-\varsigma N - (x+\mu)\cdot\nu(x) }{2(t+1)} \Big) U.
\end{displaymath}
Thanks to \refeq{mu_bom} we obtain $\mathcal{B} _\sigma (U) \geq 0$ in $N_y \times [0,\infty)$. And we have $U\geq 0$ on $\bom \setminus N_y \times (0,\infty)$. An initial data $\varphi$ with $\varphi \leq U(\cdot, 0)$ in $\adom$ permits to conclude. \\
\fin

\noindent In the case of dimension one, we use the fact that $\Omega$ is not connected. Let us write $\Omega= \RR \setminus [a,b]$ with $a<b$ in $\RR$, and let $V$ be the function defined in $\adom \times [0,\infty)$ by:
\begin{displaymath}
V(x,t)=\left\{
   						 \begin{array}{ll}
        				A (t+1)^{-\gamma} \exp{\frac{-\parallel x + \mu_1 \parallel_2^2}{4(t+1)}} & \textrm{ if } x\leq a \\
        				A (t+1)^{-\gamma} \exp{\frac{-\parallel x + \mu_2 \parallel_2^2}{4(t+1)}} & \textrm{ if } x\geq b \ ,
    						\end{array}
			  \right.
\end{displaymath}
with $A$ and $\gamma$ like in Corollary \ref{sol_mixte}, $\mu_1$ and $\mu_2$ in $\RR$ such that
\begin{displaymath}
- (a + \mu_1)  - \varsigma \geq 0 \ ,
\end{displaymath}
and
\begin{displaymath}
(b+\mu_2)   - \varsigma \geq 0 \ .
\end{displaymath}
As $\nu(a)=1$ and $\nu(b)=-1$, we obtain with \refeq{sborne} 
\begin{displaymath}
\sigma \dt V +\dn V \geq 0 \textrm{ on } ( \{a\} \cup \{b\} ) \times [0,\infty) \ .
\end{displaymath}
Following the proof of Corollary \ref{sol_mixte}, we obtain this result:

\begin{theo} Under conditions  \refeq{dissip},  \refeq{sreg},  \refeq{i-data}, \refeq{sborne}, $N=1$ and $p>3$, Problem \refeq{Ext_pbm} admits global positive solutions vanishing at infinity, for sufficiently small initial data $\varphi$.
\end{theo}

\section*{Acknowledgment}

The author would like to thank Professor Joachim von Below for helpful discussions and advices.

\bibliography{Rault-Fujita_phenomenon}

\end{document}